\input amstex.tex
\documentstyle{amsppt}
%
%
%
%
\magnification\magstephalf

\def\La{\Lambda}
\def\la{\lambda}
\def\Uq{U_{q}(\widehat{\frak{sl}}(2|1))} 
\def\UMN{\widehat{\frak{sl}}(M+1|N+1)}
\def\UqMN{U_{q}({\widehat{\frak{sl}}(M+1|N+1))}}

\def\ch{\operatorname{ch}}
\def\tr{\operatorname{tr}}

\def\ignore#1{ } 
\NoBlackBoxes
\TagsOnRight
\topmatter
\title
A level-one representation of the quantum affine superalgebra
$\UqMN$
\endtitle
\author
$^1$Kazuhiro Kimura, $^2$Jun'ichi Shiraishi and $^3$Jun Uchiyama
\endauthor
\affil
$^1$Department of Physics, Osaka Institute of Technology, Omiya,
Osaka 535, Japan \\
$^2$Institute for Solid State Physics, 
University of Tokyo, 
Roppongi, Tokyo 106, Japan\\
$^3$Department of Physics, Rikkyo University,
Nishi-Ikebukuro, Tokyo 171, Japan 
\endaffil
\abstract
A level-one representation of the quantum affine superalgebra
$\UqMN$ and vertex operators associated with the 
fundamental representations are 
constructed in terms of free bosonic fields.
Character formulas of level-one irreducible highest weight modules of 
$\Uq$ are conjectured.
\endabstract
\endtopmatter
\document
\head 1. Introduction \endhead 
In the recent studies of the mathematical physics 
and solvable systems \cite{DFJMN} (see also references in \cite{JM}), 
affine quantum algebras \cite{D1,J} played an important role.
In these works, 
techniques of the representation theory, 
infinite dimensional highest weight modules and 
corresponding vertex operators etc. have proved to be very powerful to 
study low-dimensional systems. Therefore, it may be natural to expect 
that representation theories of the affine Lie {\it superalgebras} or 
their quantum analogues will greatly help our 
future studies (as for applications to the number theory see \cite{KWk1}).
The defining relations of the quantum affine superalgebras
are obtained by Yamane \cite{Y}. 
In that paper, the Drinfeld realization is also
studied \cite{D2}. 
The representation theories of the Lie superalgebras are much more complicated
than non-super cases and have reach structures \cite{K1,K2,K3,FSS,KWn,KWk1}. 
Hence, to obtain concrete representation spaces
is desirable.
The aim of this article is to construct a
level-one representation of $\UqMN$ by bosonizing the 
Drinfeld generators basing on the free boson representation of 
$\UMN$ obtained in \cite{BCMN} and 
study character formulas for highest weight modules.

This paper is organized as follows.
In this section we review the definition of the quantum affine
superalgebra $\UqMN$ and 
construct level-zero representations. 
In Section 2, we study the bosonization of $\UqMN$ and
using that, character formulas for level-one irreducible 
highest weight modules of $\UqMN$ are conjectured.
Section 3 is devoted 
to the study of the bosonization of the vertex operators.

\subhead 1.1 
Quantum affine super algebra $\UqMN$ 
\endsubhead 
We will study the quantum affine superalgebra $\UqMN$ for
$M,N=0,1,\cdots$ and we will restrict our analysis to $M \ne N$.
The Cartan matrix of the affine Lie superalgebra $\UMN$ is 
$$
(a_{ij})=\left( \matrix
0 & -1 &  &  & &&&&&1 \\
-1 & 2 & -1 &  &  &&&&& \\
 & -1 & 2 & \ddots &  &&&&&  \\
 &  & \ddots & \ddots & \ -1 & \\ 
& & & -1 & 2 & -1 \\
& & & & -1 & 0 & 1 &  \\
& & & & & 1 & -2 &\ddots \\
& & & & & & \ddots &\ddots &1 &  \\ & & & & & & & 1 & -2 & 1 \\
1 &  &  &&&&&  & 1 & -2
\endmatrix \right) \quad (0 \leq i,j \leq M+N+1), \tag 1.1
$$
where the diagonal part is
$
(a_{ii})=(0,\overbrace{2,\cdots,2}^{M},0,\overbrace{-2,\cdots,-2}^N).
$
Let us introduce orthonormal basis 
$\{\varepsilon'_i|i=1,\cdots,M+N+2\}$ with the 
bilinear form 
$(\varepsilon'_i|\varepsilon'_j)=\nu_i \delta_{i,j}$, where 
$\nu_i=1$ for $i=1,\cdots,M+1$ and 
$\nu_j=-1$ for $i=M+2,\cdots,M+N+2$.
Define 
$
\varepsilon_i= \varepsilon'_i-
  \nu_i  \sum_{j=1}^{M+N+2}\varepsilon'_j/(M-N)
$.
The classical simple roots are 
defined by
$
\bar{\alpha}_i=\nu_i \varepsilon'_i-\nu_{i+1} \varepsilon'_{i+1}
$
and the classical weights are 
$
\bar{\Lambda}_i=\sum_{j=1}^{i} \varepsilon_j 
$
for $i=1,\cdots,M+N+1$.
Introduce the 
affine weight
$\Lambda_0$ and the null root $\delta$ having 
$(\Lambda_0|\varepsilon'_i)=(\delta|\varepsilon'_i)=0$
for $i=1,\cdots,M+N+2$ and
$(\Lambda_0|\Lambda_0)=(\delta|\delta)=0,
(\Lambda_0|\delta)=1$.
The other affine weights and affine roots are given by 
$
\Lambda_i=\bar{\Lambda}_i+\Lambda_0
$
and 
$
\alpha_i=\bar{\alpha}_i
$
for $i=1,\cdots,M+N+1$ and $\alpha_0=\delta-\sum_{i=1}^{M+N+1}\alpha_i$.
Let $P=\oplus_{i=0}^{M+N+1}\Bbb{Z}\La_i\oplus\Bbb{Z}\delta$ 
and $P^*=\oplus_{i=0}^{M+N+1} \Bbb{Z}h_i\oplus\Bbb{Z}d$ 
be the affine weight lattice and its dual lattice, respectively.

The quantum affine algebra $\UqMN$ \cite{Y} ia a $q$-analogue of the universal 
enveloping algebra of $\UMN$ generated by the Chevalley 
generators, $e_i,f_i,t_i^{\pm 1}$ over the 
base field ${\Bbb Q}(q)$. The $\Bbb{Z}_2$-grading
$|\cdot|:\UqMN\rightarrow\Bbb{Z}_2$
of the generators 
are:
$|e_0|=|f_0|=|e_{M+1}|=|f_{M+1}|=1$ and zero otherwise. 
The relations among these generators are 
$$
\align
t_it_j & =t_jt_i, \tag 1.2 \\ 
t_ie_jt_i^{-1} & =q^{a_{ij}}e_j, \tag1.3 \\ 
t_if_jt_i^{-1} &=q^{-a_{ij}}f_j, \tag1.4 \\ 
[e_i,f_j] &=\delta_{ij}{ t_i-t_i^{-1} \over q-q^{-1}}, \tag1.5
\endalign
$$
$$
\left.
\aligned
[e_j,[e_j,e_i]_{q^{-1}}]_q&=0\\
[f_j,[f_j,f_i]_{q^{-1}}]_q&=0  
\endaligned\;\;\right\}
\qquad\text{for}\;|a_{ij}| = 1,i\not=0,M+1 \tag 1.6 
$$
$$
\left.
\aligned
&[e_l ,[e_k,[e_l,e_m]_{q^{-1}}]_q]=0\\
&[f_l ,[f_k,[f_l,f_m]_{q^{-1}}]_q]=0
\endaligned\;\;\right\}
\qquad 
\aligned \text{for}\;
&(k,l,m)\\=&(M+2,M+1,M),(1,0,M+N+1),\endaligned \tag1.7
$$
where we have used the notations
$[X,Y]_\xi=XY-(-1)^{|X||Y|}\xi YX$.
Here and hereafter, we write $[X,Y]_1$ as $[X,Y]$ for simplicity.
If $M=0$ or $N=0$, we have extra fifth order Serre relations.
As for the explicit forms of the extra Serre relations,
we will refer the reader to the ref. \cite{Y}.
 
The quantum affine super algebra $\UqMN$ 
can be endowed with the graded Hopf algebra structure.
We take the following coproduct 
$$
\align
\Delta(e_i)=e_i\otimes 1+t_i\otimes e_i, \;\;
\Delta(f_i)=f_i\otimes t_i^{-1}+1\otimes f_i ,\;\;
\Delta(t_i^{\pm 1})=t_i^{\pm 1}\otimes t_i^{\pm 1},  \tag1.8 \endalign 
$$
and the antipode 
$$
a(e_i)=-t_i^{-1} e_i,\quad a(f_i)=-f_i t_i,\quad 
a(t_i^{\pm 1})=t_i^{\mp 1}. \tag1.9
$$
The coproduct is an algebra automorphism 
$\Delta(xy)=\Delta(x)\Delta(y)$ and the antipode is a graded
algebra anti-autoorphism 
$a(xy)=(-1)^{|x||y|}a(y)a(x)$ for $x,y \in \UqMN$.
Let $V$ and $W$ be graded representations of $\UqMN$.
Hereafter, the ${\Bbb Z}_2$-grading on representation space will 
also be denoted by $|\cdot|$.
The graded action of $\UqMN\otimes \UqMN$ on the tensor 
representation
$V\otimes W$
is defined by
$x\otimes y \cdot v\otimes w=(-1)^{|y||v|} xv\otimes yw$
for $x,y \in \UqMN$ and $v\in V,w\in W$

In order to construct the bosonic representation of $\UqMN$, we 
give another realization of $\UqMN$ using the Drinfeld basis \cite{Y}:
$\{X^{i,\pm}_m,h^i_n,(K^i)^{\pm 1},\gamma^{\pm 1/2}|$ $
 i=1,\cdots,M+N+1,m \in 
\Bbb{Z}, n \in \Bbb{Z}_{\ne 0}\}$ .
The $\Bbb{Z}_2$-grading
of the Drinfeld generators are:
$|X^{\pm,M+1}_m|=1$ for $m\in {\Bbb Z}$ and zero otherwise. 
The relations are
$$
\align
\gamma &\text{\;is\; central }, \tag1.10 \\ 
[K^i,h^j_m]&=0, \;\;[h^i_m,h^j_n] = \delta_{m+n,0} 
{[a_{ij}m](\gamma^m-\gamma^{-m}) \over m(q-q^{-1})}, \tag1.11 \\
K^iX^{\pm,j}_m &=q^{\pm a_{ij}}X^{\pm,j}_m K^i,\tag1.12 \\ 
[h^i_m,X^{\pm,j}_n]&=\pm {[a_{ij}m] \over m}
\gamma^{\pm |m|/2}
X^{\pm,j}_{n+m},\tag1.13 \\ 
[X^{+,i}_m,X^{-,j}_n]&=\delta_{i,j}{1 \over q-q^{-1}}
(\gamma^{(m-n)/2}\psi^{+,j}_{m+n} -\gamma^{-(m-n)/2 }
\psi^{-,j}_{m+n}),\tag1.14 \endalign 
$$
$$
\align
&[X^{\pm ,i}_m,X^{\pm ,j}_n]=0 \quad 
\text{for} \quad a_{ij}=0, \tag1.15 \\
&[X^{\pm,i}_{m+1},X^{\pm,j}_{n}]_{q^{\pm a_{ij}}} 
+[X^{\pm,j}_{n+1},X^{\pm,i}_{m}]_{q^{\pm a_{ij}}}=0, 
\quad\text{for} \quad a_{ij}\ne 0
\tag1.16 \\ 
&\text{Sym}_{l,m}[X^{\pm,i}_l,[X^{\pm,i}_m,X^{\pm, j}_n]_{q^{-1}}]_{q} 
 =0 , \qquad \text{for} \quad a_{ij}=0, i\ne M+1\tag1.17 \\ 
&\text{Sym}_{k,m} 
[X^{\pm,M+1}_k,[X^{\pm,M+2}_l,[X^{\pm,M+1}_m,X^{\pm,M}_n]_{q^{- 1}}
]_q] =0,  \tag1.18 \endalign
$$
where 
$$
\sum_{m \in \Bbb{Z}} \psi^{\pm,i}_m z^{-m} 
=(K^i)^{\pm 1}
\exp\left(\pm (q-q^{-1})\sum_{m>0}h^i_{\pm n}z^{\mp n}\right), \tag1.19 
$$
and the symbol $\text{Sym}_{k,l}$ means symmetrization with respect to
$k$ and $l$. We used the standard notation
$
[x]={q^x-q^{-x}\over q-q^{-1}}.
$
If $M=0$ or $N=0$, we have extra fifth order Serre relations for 
the Drinfeld generators.
As for the explicit forms, we will refer the reader to the ref. \cite{Y}.

The Chevalley generators are obtained by the formulas:
$$
\align
t_i&= K_i,\quad
e_i = X^{+,i}_0,\quad
f_i = X^{-,i}_0 \quad \text{for}\;i=1,\cdots,M+N+1 \tag1.20 \\ 
t_0&= \gamma K^{-1}_1 \cdots K^{-1}_{M+N+1} \tag1.21 \\ 
e_0&= (-1)^{N+1}
[X_0^{-,M+N+1}\cdots,[X_0^{-,M+2},[X_0^{-,M+1}\cdots,[X^{-,2}_0,
X^{-,1}_1]_{q{-1}}\cdots ]_{q{-1}}]_q\cdots]_q  \\
&\qquad \qquad \times (K^1\cdots K^{M+N+1})^{-1} \tag1.22 \\ 
f_0&= K_1\cdots 
K_{M+N+1}
[\cdots[[\cdots [
X^{+,1}_{-1},X^{+,2}_0]_q,\cdots X^{+,M+1}_0
]_q,X^{+,M+2}_0]_{q^{-1}},\cdots X^{+,M+N+1}_0]_{q^{-1}}
\tag1.23 \endalign
$$

\subhead 1.2 level-zero representations of $\UqMN$\endsubhead
Let $E_{i,j}$ be the $(M+N+2)\times (M+N+2)$ matrix whose $(i,j)$-element is 
unity and zero elsewhere, set
$v_i=^t (0,\cdots0,\mathop{1}\limits^{i },0,\cdots,0)$ for 
$i=1,\cdots,M+N+2$.
We will adopt the ${\Bbb Z}_2$-grading to the basis by
$|v_i|=(\nu_i+1)/2$. For the sake of simplicity, we 
will not study another possibility $|v_i|=(-\nu_i+1)/2$ in this article.
The $M+N+2$ dimensional level-zero 
representation $V_z$ of $\UqMN$ with basis 
$\{v_i\otimes z^n|i,1,\cdots,M+N+2,n\in {\Bbb Z}\}$
is defined by
$$
\align
e_i&=E_{i,i+1}, \quad
f_i=\nu_i E_{i+1,i},   \quad
t_i=q^{\nu_i E_{i,i}-\nu_{i+1}E_{i+1,i+1}}, \\
e_0&=-z E_{M+N+2,1},  \quad
f_0=z^{-1} E_{1,M+N+2},  \quad
t_0=q^{-E_{1,1}-E_{M+N+2,M+N+2}}, \tag1.24
\endalign
$$
for $i=1,\cdots,M+N+1$.
Let $V_z^*$ be the dual space of $V_z$ with basis
$\{v_i^*\otimes z^n|i,1,\cdots,M+N+2,n\in {\Bbb Z}\}$
such that 
$\langle  v_i^* \otimes z^m,v_j \otimes z^n\rangle=
\delta_{i,j}\delta_{m+n,0}$.
We also regard $v^*_i$ as the vector
$v_i=^t (0,\cdots0,\mathop{1}\limits^{i },0,\cdots,0)$.
The $\UqMN$-module structure is given by 
$\langle xv,w\rangle=\langle v,(-1)^{|x||v|}a(x)w\rangle$ for 
$v\in V_z^*, w\in V_z$ and we call the module
as $V_z^{*a}$.
The representation is:
$$
\align
e_i&=- \nu_i\nu_{i+1}q^{-\nu_i} E_{i+1,i}, \quad
f_i=-\nu_i q^{\nu_i} E_{i,i+1},   \quad
t_i=q^{-\nu_i E_{i,i}+\nu_{i+1}E_{i+1,i+1}},\\
e_0&=q z E_{1,M+N+2},  \quad
f_0=q^{-1}z^{-1} E_{M+N+2,1},  \quad
t_0=q^{E_{1,1}+E_{M+N+2,M+N+2}}. \tag1.25
\endalign
$$
The Drinfeld generators on $V_z$ are represented by
%
$$\align
h^{i}_{m} &= {[m]\over m} (q^{\mu_{i}}z)^{m}
\left(\nu_{i}q^{-\nu_{i}m}E_{i,i} -
\nu_{i+1}q^{\nu_{i+1}m}E_{i+1,i+1}\right), 
\quad   K^i=q^{\nu_{i}E_{i,i} -\nu_{i+1}E_{i+1,i+1} },\\
X^{+,i}_{m} &= (q^{\mu_{i}}z)^{m}E_{i,i+1}, \quad
X^{-,i}_{m} = \nu_{i}(q^{\mu_{i}}z)^{m}E_{i+1,i}, \tag1.26
\endalign
$$
and on $V_z^{*a}$
$$\align
h^{i}_{m} &= -{[m]\over m}
(q^{-\mu_{i}}z)^{m}\left(\nu_{i}q^{\nu_{i}m}
E_{i,i} - \nu_{i+1}q^{-\nu_{i+1}m}E_{i+1,i+1}\right), 
\quad  K^i=q^{-\nu_{i}E_{i,i} +\nu_{i+1}E_{i+1,i+1} }   \\
X^{+,i}_{m} &= -\nu_{i}\nu_{i+1}q^{-\nu_{i}}(q^{-\mu_{i}}z)^{m}E_{i+1,i}, 
\quad
X^{-,i}_{m} = -\nu_{i}q^{\nu_{i}}(q^{-\mu_{i}}z)^{m}E_{i,i+1}, \tag1.27
\endalign
$$
where $\mu_{i} = \sum_{k=1}^{i} \nu_{k}$.

\head 2. A level-one representation of $\UqMN$\endhead 
\subhead 2.1 Free boson realization \endsubhead
Now, we will study the free boson realization of $\UqMN$
which gives us a level-one representation. 
It is well known that the representations of 
the non-super affine algebras are constructed in terms of 
bosonic fields at level-one \cite{FK}. 
Basing on this realization, Frenkel and Jing constructed the 
free boson representation of the quantum affine algebras 
at level-one \cite{FJ}. 
We will show that 
this kind of bosonization 
can be extended to the affine superalgebras of $A$-type. 
Our representation can be regarded as a $q$-deformation of the 
free field realization of $\UMN$ studied by Bouwknegt et.al. \cite{BCMN}.
The structure of the deformation is 
essentially the same as that of Frenkel and Jing's except for the 
deformation of the $\beta$-$\gamma$ ghost-system.
To this ghost-system, however, 
the technique for bosonizing 
the $\beta$-$\gamma$ ghost
which was discussed in the papers on 
deformed Wakimoto realization of 
$U_{q}(\widehat{\frak{sl}}(N))$ 
(see \cite{AOS} and references therein) is applicable. 

Let us introduce the bosonic oscillators
$
\{a^i_n,b^j_n,c^j_n,Q_{a^i},Q_{b^j},Q_{c^j}|
n \in {\Bbb Z}, i=1,\cdots,M+1,j=1,\cdots,N+1\}$
satisfying the commutation relations
$$
\align
[a^i_m,a^j_n]&=\delta_{i,j}\delta_{m+n,0}{[m]^2 \over m},\qquad 
[a^i_0,Q_{a^j}] = \delta_{i,j},\tag2.1 \\
[b^i_m,b^j_n]&=-\delta_{i,j}\delta_{m+n,0}{[m]^2 \over m},\qquad 
[b^i_0,Q_{b^j}] = -\delta_{i,j} ,\tag2.2 \\
[c^i_m,c^j_n]&=\delta_{i,j}\delta_{m+n,0}{[m]^2 \over m},\qquad 
[c^i_0,Q_{c^j}] = \delta_{i,j}. \tag2.3
\endalign
$$
The remaining commutators vanish. 

Define the generating functions for the Drinfeld basis by
$X^{\pm,i}(z)=\sum_{m\in {\Bbb Z}}X^{\pm,i}_mz^{-m-1} $,
and introduce $h^i_0$ by setting $K^i=q^{h^i_0}$.
Define
$
Q_{h^i}=Q_{a^i}-Q_{a^{i+1}}
$ for $i=1,\cdots,M$,
$
Q_{h^{M+1}}=Q_{a^{M+1}}+Q_{b^{1}}
$
and
$
Q_{h^{M+1+j}}=-Q_{b^j}+Q_{b^{j+1}}
$ for $j=1,\cdots,N$.
Let us introduce the notation
$$
h^i(z;\beta)=-\sum_{n \ne 0}{h^i_n \over 
[n]}q^{-\beta|n|}z^{-n} +Q_{h^i}+h^i_0\ln z, \tag2.4 
$$
for 
the Drinfeld generators $h^i_m$, $Q_{h^i}$ and $\beta\in {\Bbb R}$.
In this article, 
we will adopt this notation for other bosonic fields, for example,
the boson field $c^j(z;\beta)$ should be defined in the same way.
We introduce the $q$-differential operator defined by $
_1\partial_z f(z)={f(qz)-f(q^{-1}z) \over (q-q^{-1})z}$.
Now we state the result of the bosonization.
\proclaim{Proposition 2.1}
The Drinfeld generators at level-one are realized by the free boson fields as
$$
\align
\gamma &=\; q,\tag2.5 \\
h^i_m &= \; a^i_mq^{-|m|/2}-a^{i+1}_mq^{|m|/2},
\tag2.6 \\ 
h^{M+1}_m &=  \; a^{M+1}_mq^{-|m|/2}+b^1_mq^{-|m|/2}, 
\tag2.7 \\ 
h^{M+1+j}_m &=  \; -b^j_mq^{|m|/2}+b^{j+1}_mq^{-|m|/2},
\tag2.8 \\
X^{+,i}(z) &=  \; :e^{h^i(z;1/2)}:e^{i\pi a^i_0},
\tag2.9 \\ 
X^{+,M+1}(z) &=  \; :e^{h^{M+1}(z;1/2)}e^{c^1(z;0)}: 
\prod^M_{i=1}e^{-i\pi a^i_0}, 
\tag2.10 \\ 
X^{+,M+1+j}(z) &=  \; :e^{h^{M+1+j}(z;1/2)}[ _1\partial_z e^{-c^j(z;0)}]
e^{c^{j+1}(z;0)}: ,
\tag2.11 \\ 
X^{-,i}(z) &=  \; -:e^{-h^i(z;-1/2)}:e^{-i\pi a^i_0}, 
\tag2.12 \\ 
X^{-,M+1}(z) &=  \; :e^{-h^{M+1}(z;-1/2)}[_1\partial_ze^{-c^1(z;0)}]:  
\prod^M_{i=1}e^{i\pi a^i_0}, 
\tag2.13 \\ 
X^{-,M+1+j}(z) &=  \;-:e^{-h^{M+1+j}(z;-1/2)}e^{c^{j}(z;0)}[ 
_1\partial_z e^{-c^{j+1}(z;0)}]:, 
\tag2.14 \endalign
$$
for $m\in {\Bbb Z}_{\ne 0}$, $i=1,\cdots,M$ and $j=1,\cdots,N$.
The usual normal ordering is denoted by $:\cdots:$.
\endproclaim
We can check the commutation relations by studying the operator
products among the bosonized generators.

\subhead 2.2 Highest weight $\Uq$-modules \endsubhead
To exploit special features of the quantum affine superalgebras,
we study the simplest example $\Uq$ of the level-one representation
obtained in the last subsection.
We begin by defining the Fock module.  
The vacuum vector $ |0 \rangle$ is defined by
$
a^{i}_{n} |0 \rangle = b_{n} |0 \rangle = 
 c_{n} |0 \rangle = 0
$
for
$
n\ge 0, 
$
and the vector carrying the weight 
$(\la_{a^{1}},\la_{a^{2}},\la_{b},\la_{c}) \in \Bbb{C}^{4}$ by
$$
|\la_{a^{1}},\la_{a^{2}},\la_{b},\la_{c}\rangle = e^{\la_{a^{1}}Q_{a^{1}}+
\la_{a^{2}}Q_{a^{2}}+\la_{b}Q_{b}+\la_{c}Q_{c}}|0\rangle. \tag 2.15
$$
The Fock module ${\Cal F}_{\la_{a^{1}},\la_{a^{2}},\la_{b},\la_{c}}$ 
is generated by acting creation operators 
$h^1_n=a_{n}^{1}q^{n/2}-a^2_n q^{-n/2},
 h^2_n=a_{n}^{2}q^{n/2}+b_n   q^{n/2}$ and 
$c_{n}\; (n<0)$ on $|\la_{a^{1}},\la_{a^{2}},\la_{b},\la_{c}\rangle$.

To obtain highest weight vectors of $\Uq$, we impose the conditions:
$$
\aligned 
 &e_i|\la_{a^{1}},\la_{a^{2}},\la_{b},\la_{c}\rangle = 0,\\ 
&h_i|\la_{a^{1}},\la_{a^{2}},\la_{b},\la_{c}\rangle = 
\lambda^i|\la_{a^{1}},\la_{a^{2}},\la_{b},\la_{c}\rangle,
\endaligned
\qquad \text{for}\;i=0,1,2, \tag2.16
$$
Solving these equations, we obtain the following classification:
\roster
\item 
$(\la_{a^{1}},\la_{a^{2}},\la_{b},\la_{c})=
(\beta,\beta,\beta-\alpha,-\alpha)$, where $\alpha$ and 
$\beta$ are arbitrary. The weight of this vector is  
$(\lambda^0,\lambda^1,\lambda^2)=(1-\alpha,0,\alpha)$.
Thus we have the identification:
$|(1-\alpha)\Lambda_0+\alpha \Lambda_2\rangle
=|\beta,\beta,\beta-\alpha,-\alpha\rangle$.
\item $(\la_{a^{1}},\la_{a^{2}},\la_{b},\la_{c})=
(\beta+1,\beta,\beta,0)$, where $\beta$ is arbitrary. 
The weight is  
$(\lambda^0,\lambda^1,\lambda^2)=(0,1,0)$. We have 
$|\Lambda_1\rangle=|\beta+1,\beta,\beta,0\rangle$.
\item $(\la_{a^{1}},\la_{a^{2}},\la_{b},\la_{c})=
(\beta+1,\beta+1,\beta,0)$, where $\beta$ is arbitrary. 
The weight is  
$(\lambda^0,\lambda^1,\lambda^2)=(0,0,1)$ and
we have 
$|\Lambda_2\rangle=|\beta+1,\beta+1,\beta,0\rangle$.
\endroster
According to this classification, let us introduce the following spaces
$$
\align
&{\Cal F}_{(\alpha;\beta)}=
\bigoplus_{i,j\in {\Bbb Z}} 
{\Cal F}_{\beta+i,\beta-i+j,\beta-\alpha+j,-\alpha+j}, \tag2.17 \\
{\Cal F}_{((1,0);\beta)}=
\bigoplus_{i,j\in {\Bbb Z}} &
{\Cal F}_{\beta+1+i,\beta-i+j,\beta+j,j},\quad
{\Cal F}_{((0,1);\beta)}=
\bigoplus_{i,j\in {\Bbb Z}} 
{\Cal F}_{\beta+1+i,\beta+1-i+j,\beta+j,j}.
\endalign
$$
It is not difficult to see that the bosonized actions of $\Uq$ on these spaces
are closed i.e. 
$\Uq {\Cal F}_{(*;\beta)}={\Cal F}_{(*;\beta)}$ where 
$*=\alpha,(1,0),(0,1)$. 
These spaces are not irreducible in general. 
It is convenient to introduce a pair of Fermionic fields 
$\eta(z)=
 \sum_{n\in {\Bbb Z}}\eta_{n} z^{-n-1} = :e^{c(z;0)}:$
and 
$\xi(z)
= \sum_{n\in {\Bbb Z }}\xi_{n} z^{-n}= :e^{-c(z)}:
$ to obtain the irreducible subspaces of these
$\Uq$-modules.
The mode expansion of $\eta(z)$, $\xi(z)$ is well defined 
on ${\Cal F}_{(\alpha;\beta)}$ for $\alpha\in {\Bbb Z}$ and 
on ${\Cal F}_{((1,0);\beta)}$, ${\Cal F}_{((0,1);\beta)}$ and the
relations are 
$
\{ \xi_{r} , \xi_{s} \} = \{\eta_{r}, \eta_{s} \} = 0,
\{ \xi_{r},\eta_{s}\} = \delta_{r+s,0}. 
$
In these cases, 
we have the direct sum decompositions
${\Cal F}_{(*;\beta)}=
\eta_0\xi_0{\Cal F}_{(*;\beta)}
\oplus
\xi_0\eta_0{\Cal F}_{(*;\beta)}.
$
As usual, we call $\eta_0\xi_0{\Cal F}_{(*;\beta)}$
as $\text{Ker}\,{\eta_0}_{(*;\beta)}$ 
and 
${\Cal F}_{(*;\beta)}/\eta_0\xi_0{\Cal F}_{(*;\beta)}$
as $\text{Coker}{\,\eta_0}_{(*;\beta)}$.
Since $\eta_0$ commutes (or anti-commutes) with
every element of $\Uq$, we can regard 
$\text{Ker}\,{\eta_0}_{(*;\beta)}$ and 
$\text{Coker}\,{\eta_0}_{(*;\beta)}$ 
as $\Uq$-modules. 

From now on, we study the character formulas of these $\Uq$-modules we have
constructed in the bosonic Fock space.
The character of a space ${\Cal F}$ is defined by 
$$
\ch_{\Cal F}(q,x,y) \equiv \tr_{\Cal F} q^{-d}x^{h^{1}_0}y^{h^{2}_0}, \tag2.18
$$
where the grading operator $d$ is \cite{BCMN}
$$\align
d &=- \sum_{n\ge 1} (\frac{n}{[n]})^{2}\{a^{1}_{-n}a^{1}_{n} 
+a^{2}_{-n}a^{2}_{n} - b_{-n}b_{n} + c_{-n}c_{n} - (a^{1}_{-n} 
+a^{2}_{-n} +b_{-n} )
(a^{1}_{n} +a^{2}_{n} +b_{n} ) \} \\ 
&-\frac{1}{2}\{(a^{1}_{0})^{2} + (a^{2}_{0})^{2} -(b_{0})^{2} 
+c_{0}(c_{0}+1) -(a^{1}_{0} +a^{2}_{0} +b_{0})^{2}  \},\tag2.19  \endalign 
$$
and $h^1_0=a^1_0-a^2_0$, $h^2_0=a^2_0+b_0$.

(I) {\it Character of ${\Cal F}_{(\alpha;\beta)}$ for $\alpha\not\in 
{\Bbb Z}$.}\quad Since $\eta_0$ is not defined on this module,
it is expected that ${\Cal F}_{(\alpha;\beta)}$ is irreducible 
highest weight $\Uq$-module. Thus, we conjecture the 
following.
\proclaim{Conjecture 2.1}
We have the identification of the highest weight $\Uq$-modules:
$$
{\Cal F}_{(\alpha;\beta)}\cong 
V((1-\alpha)\Lambda_0+\alpha \Lambda_2)\quad \text{for}\; 
\alpha\not\in {\Bbb Z}\;\text{and\;arbitrary }\beta ,\tag2.20
$$
where $V(\lambda)$ denotes the 
irreducible highest weight $\Uq$-modules with the highest weight 
weight $\lambda$.
\endproclaim
\proclaim{Proposition 2.2}
We obtain the character of ${\Cal F}_{(\alpha;\beta)}$ as
$$\align
\text{ch}_{{\Cal F}_{(\alpha;\beta)} }(q,x,y)
&= 
 \frac{q^{-\frac{1}{2} \alpha (\alpha +1)}}
{\prod_{n=1}^{\infty} (1-q^{n})^{3}}\sum_{s,t\in \Bbb{Z}}
q^{\frac{1}{2}(2s^{2}-2st +t^{2} +t )}
x^{2s-t}y^{\alpha -s}.\tag2.21 
\endalign 
$$
\endproclaim

(II) {\it Character of $\text{Ker}\,{\eta_0}$ and 
$\text{Coker}\,{\eta_0}$.}\quad 
Next, let us consider the character of the modules 
on which $\eta_0$ is well defined.
\proclaim{Proposition 2.3}
The character of $\text{Ker}\,{\eta_0}_{(\alpha;\beta)}$ 
for $\alpha \in {\Bbb Z}$ is obtained as
$$\align
&\text{ch}_{\text{Ker}\,{\eta_0}_{(\alpha;\beta)}}(q,x,y)\\
&\quad =
\frac{q^{-\frac{1}{2} \alpha (\alpha +1)}}
{\prod_{n=1}^{\infty} (1-q^{n})^{3}}
\bigg(
\sum_{{s,t,l \in {\Bbb Z}}\atop  t<\alpha,l \ge 1}
(-1)^{l+1} q^{\frac{l(l-1)}{2}+l(\alpha-t)
+\frac{1}{2}(2s^{2}-2st+t^{2}+t )}
x^{2s-t}y^{\alpha -s} \\
&
\qquad\qquad\qquad\qquad+
\sum_{{s,t,l \in {\Bbb Z}} \atop t\ge \alpha ,l \ge 0}
(-1)^{l} q^{\frac{l(l+1)}{2}-l(\alpha-t)
+\frac{1}{2}(2s^{2}-2st+t^{2}+t )}
x^{2s-t}y^{\alpha -s} \bigg). \tag2.22
\endalign 
$$
The character of $\text{Coker}\,{\eta_0}_{(\alpha;\beta)}$ 
($\alpha \in {\Bbb Z}$) is 
$$\align
&\text{ch}_{\text{Coker}\,{\eta_0}_{(\alpha;\beta)}}(q,x,y)\\
&\quad =
\frac{q^{-\frac{1}{2} \alpha (\alpha +1)}}
{\prod_{n=1}^{\infty} (1-q^{n})^{3}}
\bigg(
\sum_{{s,t,l \in {\Bbb Z}}\atop t<\alpha ,l \ge 0}
(-1)^l q^{\frac{l(l-1)}{2}+l(\alpha-t)
+\frac{1}{2}(2s^{2}-2st+t^{2}+t )}
x^{2s-t}y^{\alpha -s} \\
&
\qquad\qquad\qquad\qquad+
\sum_{{s,t,l \in {\Bbb Z}}\atop t\ge \alpha ,l \ge 1}
(-1)^{l+1} q^{\frac{l(l+1)}{2}-l(\alpha-t)
+\frac{1}{2}(2s^{2}-2st+t^{2}+t )}
x^{2s-t}y^{\alpha -s} \bigg). \tag2.23
\endalign 
$$
\endproclaim
These formulas are obtained by inserting
the projectors $\eta_0\xi_0$ and $\xi_0\eta_0$ 
to the trace of the Fock space. As for the details of 
this technique, see \cite{BCMN}.
We also have the following formulas.
\proclaim{Proposition 2.4}
We have the equalities
$
\text{ch}_{\text{Coker}\,{\eta_0}_{((0,1);\beta)}}(q,x,y)=
\text{ch}_{\text{Coker}\,{\eta_0}_{(1;\beta)}}(q,x,y),
$
and
$$\align
&\text{ch}_{\text{Coker}\,{\eta_0}_{((1,0);\beta)}}(q,x,y)\\
&\quad =
\frac{1}
{\prod_{n=1}^{\infty} (1-q^{n})^{3}}
\bigg(
\sum_{{s,t,l \in {\Bbb Z}}\atop t<0 , l \ge 0}
(-1)^l q^{\frac{l(l-1)}{2}-lt
+\frac{1}{2}(2s^{2}-2st+t^{2}+t )}
x^{1+2s-t}y^{-s} \\
&
\qquad\qquad\qquad\qquad+
\sum_{{s,t,l \in {\Bbb Z}}\atop t\ge 0 , l \ge 1}
(-1)^{l+1} q^{\frac{l(l+1)}{2}+lt
+\frac{1}{2}(2s^{2}-2st+t^{2}+t )}
x^{1+2s-t}y^{ -s} \bigg). \tag2.24
\endalign 
$$
\endproclaim
Since we have the conditions 
$\eta_0|\beta,\beta,\beta-\alpha,-\alpha\rangle \ne 0$
for $\alpha=0,1,\cdots$ and 
$\eta_0|\beta,\beta,\beta-\alpha,-\alpha\rangle =0$
for $\alpha=-1,-2,\cdots$, the modules 
$\text{Coker}\,{\eta_0}_{(\alpha;\beta)}$ ($\alpha=0,1,\cdots$),
$\text{Coker}\,{\eta_0}_{((0,1);\beta)}$,
$\text{Coker}\,{\eta_0}_{((1,0);\beta)}$ and  
$\text{Ker}\,{\eta_0}_{(\alpha;\beta)}$ ($\alpha=-1,-2,\cdots$)
 are highest weight $\Uq$-modules.
It is expected that these modules are also 
irreducible with respect to the action of $\Uq$.
\proclaim{Conjecture 2.2}
We have
the following identifications of the highest weight 
$\Uq$-modules:
$$
\align
V((1-\alpha)\Lambda_0+\alpha \Lambda_2)
&\cong 
\text{Coker}\,{\eta_0}_{(\alpha;\beta)}\quad \text{for}\; 
\alpha=0,1,\cdots,\\
&\cong 
\text{Ker}\,{\eta_0}_{(\alpha;\beta)}\quad \text{for}\; 
\alpha=-1,-2,\cdots, \tag2.25
\endalign
$$
and 
$
V(\Lambda_1)\cong 
\text{Coker}\,{\eta_0}_{((1,0);\beta)}$, 
$
V(\Lambda_2)\cong 
\text{Coker}\,{\eta_0}_{((0,1);\beta)}$ for arbitrary $\beta$.
\endproclaim

We have checked the validity of Conjecture 2.1 and 
Conjecture 2.2 for $V(\Lambda_0),V(\Lambda_1),V(\Lambda_2)$
up to certain degrees by comparing our results with
those of 
Kac and Wakimoto \cite{KWk2}.

\head 3. Vertex Operators \endhead
\subhead 3.1 Vertex operators for $\UqMN$ \endsubhead
In this section, we study free boson realization of 
vertex operators for $\UqMN$.
Let $V(\lambda)$ be the highest weight $\UqMN$-module
with the highest weight$\lambda$.
The ${\Bbb Z}_2$-gradation of $V(\lambda)$ is also denoted by $|\cdot|$.
The vertex operators 
$\Phi_{\lambda}^{\mu V}(z),\Phi_{\lambda}^{\mu V^*}(z),
\Psi_{\lambda}^{V \mu }(z),\Psi_{\lambda}^{V^* \mu}(z)$
are defined as the following intertwiners of  
$\UqMN$-modules if they exist: 
$$\align
&\Phi_{\lambda}^{\mu V}(z) :
 V(\lambda) \longrightarrow V(\mu)\otimes V_{z} , \quad
\Phi_{\lambda}^{\mu V^*}(z) :
 V(\lambda) \longrightarrow V(\mu)\otimes V_{z}^{*a} , 
\tag 3.1 \\
&\Psi_{\lambda}^{V \mu}(z) : 
 V(\lambda) \longrightarrow V_{z}\otimes V(\mu),\quad 
\Psi_{\lambda}^{V^* \mu}(z) : 
 V(\lambda) \longrightarrow V_{z}^{*a}\otimes V(\mu), \tag 3.2\\
&
\Phi_{\lambda}^{\mu V}(z)\cdot x=\Delta(x)\cdot  
\Phi_{\lambda}^{\mu V}(z),\quad
\Phi_{\lambda}^{\mu V^*}(z)\cdot x=\Delta(x)\cdot  \Phi_{\lambda}^{\mu
  V^*}(z), \tag3.3 \\
&\Psi_{\lambda}^{V \mu}(z)\cdot x=\Delta(x)\cdot  
\Psi_{\lambda}^{V \mu}(z) ,\quad
\Psi_{\lambda}^{V^* \mu}(z)\cdot x=\Delta(x)\cdot  \Psi_{\lambda}^{V^*
  \mu}(z) ,\tag3.4 
\endalign
$$
for $^\forall x\in \UqMN$ together with the gradation
$|\Phi_{\lambda}^{\mu V}(z)|=|\Phi_{\lambda}^{\mu V^*}(z)|=
|\Psi_{\lambda}^{V \mu}(z)|=|\Psi_{\lambda}^{V^* \mu }(z)|=0$.
We expand the vertex operators as
$
\Phi_{\lambda}^{\mu V}(z)
=
\sum_{l=1}^{M+N+2}\Phi_{\lambda l}^{\mu V}(z)\otimes v_l,
$
$
\Phi_{\lambda}^{\mu V^*}(z)
=
\sum_{l=1}^{M+N+2}\Phi_{\lambda l}^{\mu V^*}(z)\otimes v_l^*,
$
$
\Psi_{\lambda}^{V \mu }(z)
=
\sum_{l=1}^{M+N+2}v_l\otimes \Phi_{\lambda l}^{V \mu }(z),
$
and
$
\Psi_{\lambda}^{V^* \mu }(z)
=
\sum_{l=1}^{M+N+2}v_l^* \otimes \Psi_{\lambda l}^{V^* \mu}(z).
$
We define the graded action of these expanded operators by
$
\Phi_{\lambda}^{\mu V}(z) |u\rangle
=
\sum_{l=1}^{M+N+2}\Phi_{\lambda l}^{\mu V}(z)|u\rangle \otimes v_l
(-1)^{|v_l|\,||u\rangle|},
$
$
\Phi_{\lambda}^{\mu V^*}(z)|u\rangle 
=
\sum_{l=1}^{M+N+2}\Phi_{\lambda l}^{\mu V^*}(z)|u\rangle \otimes v_l^*
(-1)^{|v_l^*|\,||u\rangle|},
$
$
\Psi_{\lambda}^{V \mu }(z)|u\rangle
=
\sum_{l=1}^{M+N+2}v_l\otimes \Phi_{\lambda l}^{V \mu }(z)|u\rangle,
$
and
$
\Psi_{\lambda}^{V^* \mu }(z)|u\rangle
=
\sum_{l=1}^{M+N+2}v_l^* \otimes \Psi_{\lambda l}^{V^* \mu}(z)|u\rangle,
$
for $|u\rangle \in V(\lambda)$.

Let us introduce the following combinations of the Drinfeld operators
$$
h^{*i}_m=\sum_{j=1}^{M+N+1}
{[\alpha_{ij}m][\beta_{ij}m]\over[(M-N)m][m]}h^j_m, \;\;
Q^{*}_{h^i}=\sum_{j=1}^{M+N+1}
{\alpha_{ij} \beta_{ij}\over M-N}Q_{h^j},\;\;
h^{*i}_0=\sum_{j=1}^{M+N+1}
{\alpha_{ij} \beta_{ij}\over M-N} h^j, \tag3.5
$$
where $h^j_0$ is defined by $K^j=q^{h^j_0}$ and 
$$
\align
\alpha_{ij}&=
\left\{
\aligned
&\text{min}(i,j) \qquad\qquad\qquad \text{if}\;\text{min}(i,j)\leq M+1,\\
&2(M+1)-\text{min}(i,j) \qquad \text{if}\;\text{min}(i,j)> M+1,
\endaligned
\right. \tag3.6 \\
\beta_{ij}&=
\left\{
\aligned
&M-N-\text{max}(i,j) \qquad\qquad\qquad \text{if}\;\text{max}(i,j)\leq M+1,\\
&-M-N-2+\text{max}(i,j) \qquad \text{if}\;\text{max}(i,j)> M+1. 
\endaligned 
\right. \tag3.7
\endalign
$$
Note that using these notations, we have the inverse of the 
Cartan matrix as 
$
(a_{ij})^{-1}=\alpha_{ij}\beta_{ij}/(M-N).
$
We obtain the relations
$
[h^{*i}_m,h^j_n]=\delta_{i,j}\delta_{m+n,0}{[m]^2\over m},\;\;
[h^{*i}_m,h^{*j}_n]=\delta_{m+n,0}{[a^{-1}_{ij} m][m]\over m},
$
and 
$
[h^{*i},Q_{h^j}]=\delta_{i,j},
$
$
[h^{*i},Q^*_{h^{j}}]=a^{-1}_{ij}
$.

Define the operators $\phi_l(z),\phi^*_l(z),\psi_l(z)$, and  $\psi^*_l(z)$
$(i=1,\cdots,M+N+2)$ 
iteratively by 
$$
\align
\phi_{M+N+2}(z)&=\;
:e^{-h^*_{M+N+1}(q^{M-N+1}z;-1/2)+c^{N+1}(q^{M-N+1}z;0)}: 
\prod_{i=1}^{M+1}e^{i\pi {1-k \over M-N}a^i_0},\tag3.8\\
\nu_l \phi_l(z)&=[\phi_{l+1}(z),f_l]_{q^{\nu_{l+1}}},\tag3.9\\
%
\phi^*_{1}(z)&=\;
:e^{h^*_{1}(qz;-1/2)}: 
\prod_{i=1}^{M+1}e^{i\pi {k-1 \over M-N}a^i_0},\tag3.10 \\
-\nu_l q^{\nu_l} \phi^*_{l+1}(z)&=[\phi^*_{l}(z),f_l]_{q^{\nu_{l}}},
\tag3.11 \\
%
\psi_{1}(z)&=\;
:e^{-h^*_{1}(qz;1/2)}: 
\prod_{i=1}^{M+1}e^{i\pi {1-k \over M-N}a^i_0},\tag3.12 \\
\psi_{l+1}(z)&=[\psi_{l}(z),e_l]_{q^{\nu_{l}}},\tag3.13 \\
%
\psi^*_{M+N+2}(z)&=\;
:e^{h^*_{M+N+1}(q^{-M+N+1}z;1/2)}
[ _1\partial_z e^{-c^{N+1}(q^{-M+N+1}z;0)}]:  
\prod_{i=1}^{M+1}e^{i\pi {k-1 \over M-N}a^i_0},\tag3.14 \\
-\nu_l\nu_{l+1}q^{-\nu_l} \psi^*_l(z)&=[\psi^*_{l+1}(z),e_l]_{q^{\nu_{l+1}}},
\tag3.15 \endalign
$$
where $h^*_i(z;\beta)$ is defined in the same manner as
$h_i(z;\beta)$.
The gradations are given by 
$|\phi_l(z) |= |\phi^*_l(z)| = |\psi_l(z)| =| \psi^*_l(z)| =
{\nu_l+1\over 2}$.
Define the operators $\phi(z),\phi^*(z),\psi(z)$ and $\psi^*(z)$ by
$\phi(z)=\sum_{l=1}^{M+N+2}\phi_l(z)\otimes v_l$,
$\phi^*(z)=\sum_{l=1}^{M+N+2}\phi^*_l(z)\otimes v^*_l$,
$\psi(z)=\sum_{l=1}^{M+N+2}v_l \otimes \psi_l(z)$ and
$\psi^*(z)=\sum_{l=1}^{M+N+2}v^*_l \otimes\psi^*_l(z)$ respectively.
Then we have the following result.
\proclaim{Proposition 3.1}
The operators $\phi(z),\phi^*(z),\psi(z)$ and $\psi^*(z)$ satisfy the 
same commutation relations as 
$\Phi_{\lambda}^{\mu V}(z)$, 
$\Phi_{\lambda}^{\mu V^*}(z)$,
$\Psi_{\lambda}^{V \mu }(z)$ and
$\Psi_{\lambda}^{V^* \mu }(z)$ respectively have, respectively.
\endproclaim
To prove the proposition, 
the equations 
$
[[\psi_1(z),e_1]_q,e_1]_{q^{-1}}=0
$,
$
[\psi_1(z),e_i]=0
$ $(i\ne 1)$, 
$(q z-q^{-1}x)\psi_1(z)X^{+,1}(x)=(z-x)X^{+,1}(x)\psi_1(z)$,
and similar formulas for $\psi^*_{M+N+2}$, $\phi_{M+N+2}$ and
$\phi^*_{1}$
are helpful.

{\it Remark.} These operators can almost be determined by the method 
used for the level-one bosonization of the vertex operators of 
$U_q(\widehat{\frak{sl}}(N))$ \cite{JMMN,Ko}. Namely, we obtained
those by
studying the commutation relations 
between the vertex operators and some of the Drinfeld basis.
Relevant explicit coproduct formulas for the Drinfeld basis 
can be obtained in the same way as Chari and Pressley's \cite{CP}.
We have the bosonic feilds $c$'s 
whose contribution to the vertex operators cannot be 
determined by studying the commutation relations with $h^i_m$,
because they do not contain $c$'s. 
However, the following two information 
enables us to find the unique solutions as above:
{\it i)} the vertex operators of $q\rightarrow 1$ limit,
{\it ii)} the commutation relations 
with $X^{+,i}(z)$ (or $X^{-,i}(z)$) for type I (or type II).


\subhead 3.2 $\Uq$ case \endsubhead
We study the action of the bosonized vertex operators
of $\Uq$ 
on the Fock space defined in Subsection 2.2.
Using the bosonic representations of the vertex operators, we have the 
homomorphisms of $\Uq$-modules:
$$
\align
&
\phi(z):
\left\{
\aligned
&{\Cal F}_{(\alpha;\beta)}\rightarrow 
 {\Cal F}_{(\alpha-1;\beta+1)}\otimes V_z, \\
&{\Cal F}_{((1,0);\beta)}\rightarrow 
 {\Cal F}_{(0;\beta+1)}\otimes V_z,    \\
&{\Cal F}_{((0,1);\beta)}\rightarrow 
 {\Cal F}_{(1;\beta+1)}\otimes V_z,    \\
&{\Cal F}_{((0,1);\beta)}\rightarrow 
 {\Cal F}_{((1,0);\beta+1)}\otimes V_z,   \\
&{\Cal F}_{(3;\beta)}\rightarrow 
 {\Cal F}_{((0,1);\beta+1)}\otimes V_z,    \\
&{\Cal F}_{(2;\beta)}\rightarrow 
 {\Cal F}_{((1,0);\beta+1)}\otimes V_z,    
\endaligned
\right.\quad
\psi(z):
\left\{
\aligned
&{\Cal F}_{(\alpha;\beta)}\rightarrow 
 V_z \otimes {\Cal F}_{(\alpha-1;\beta+1)},    \\
& {\Cal F}_{((1,0);\beta)}\rightarrow 
 V_z \otimes {\Cal F}_{(0;\beta+1)},    \\
& {\Cal F}_{((0,1);\beta)}\rightarrow 
 V_z \otimes {\Cal F}_{(1;\beta+1)},    \\
&{\Cal F}_{((0,1);\beta)}\rightarrow 
 V_z \otimes {\Cal F}_{((1,0);\beta+1)},    \\
&{\Cal F}_{(3;\beta)}\rightarrow 
 V_z \otimes {\Cal F}_{((0,1);\beta+1)},    \\
&{\Cal F}_{(2;\beta)}\rightarrow 
 V_z \otimes {\Cal F}_{((1,0);\beta+1)},  \\
\endaligned
\right. \tag3.16
\\
&
\phi^*(z):
\left\{
\aligned
&{\Cal F}_{(\alpha;\beta)}\rightarrow 
 {\Cal F}_{(\alpha+1;\beta-1)}\otimes V^{*a}_z,    \\
&{\Cal F}_{((1,0);\beta)}\rightarrow 
 {\Cal F}_{(2;\beta-1)}\otimes V^{*a}_z,    \\
&{\Cal F}_{((0,1);\beta)}\rightarrow 
 {\Cal F}_{(3;\beta-1)}\otimes V^{*a}_z,    \\
&{\Cal F}_{((1,0);\beta)}\rightarrow 
 {\Cal F}_{((0,1);\beta-1)}\otimes V^{*a}_z,    \\
&{\Cal F}_{(1;\beta)}\rightarrow 
 {\Cal F}_{((0,1);\beta-1)}\otimes V^{*a}_z,    \\
&{\Cal F}_{(0;\beta)}\rightarrow 
 {\Cal F}_{((1,0);\beta-1)}\otimes V^{*a}_z, 
\endaligned
\right.\quad
\psi^*(z):
\left\{
\aligned
&{\Cal F}_{(\alpha;\beta)}\rightarrow 
 V^{*a}_z \otimes {\Cal F}_{(\alpha+1;\beta-1)},    \\
&{\Cal F}_{((1,0);\beta)}\rightarrow 
 V^{*a}_z \otimes {\Cal F}_{(2;\beta-1)},    \\
&{\Cal F}_{((0,1);\beta)}\rightarrow 
 V^{*a}_z \otimes {\Cal F}_{(3;\beta-1)},    \\
&{\Cal F}_{((1,0);\beta)}\rightarrow 
 V^{*a}_z \otimes {\Cal F}_{((0,1);\beta-1)},    \\
&{\Cal F}_{(1;\beta)}\rightarrow 
 V^{*a}_z \otimes {\Cal F}_{((0,1);\beta-1)},    \\
&{\Cal F}_{(0;\beta)}\rightarrow 
 V^{*a}_z \otimes {\Cal F}_{((1,0);\beta-1)},  
\endaligned
\right.\tag3.17
\endalign
$$

Next let us consider the vertex operators 
which intertwine highest weight $\Uq$-modules by using the 
above results.
It is easy to see that the vertex operators also 
commute (or anti-commute) with $\eta_0$. Noting this property,
the above homomorphisms and Conjecture 2.1, 2.2,
we can study the conditions of existence for the vertex operators
which intertwine the irreducible higest weight modules.
\proclaim{Conjecture 3.1}
The following vertex operators associated with the 
level-one irreducible highest weight modules exist:
$$
\align
&
\aligned
&\Phi_{\lambda_{\alpha}}^{\lambda_{\alpha-1}V}(z) :
 V(\lambda_{\alpha}) \longrightarrow V(\lambda_{\alpha-1})\otimes V_{z} , \\
&\Phi_{\Lambda_1}^{\Lambda_0 V}(z) :
 V(\Lambda_1) \longrightarrow V(\Lambda_0)\otimes V_{z} , \\
&\Phi_{\Lambda_2}^{\Lambda_1 V}(z) :
 V(\Lambda_2) \longrightarrow V(\Lambda_1)\otimes V_{z} , 
\endaligned
\qquad
\aligned
&\Psi_{\lambda_{\alpha}}^{V \lambda_{\alpha-1}}(z) : 
 V(\lambda_{\alpha}) \longrightarrow V_{z}\otimes V(\lambda_{\alpha-1}),\\
&\Psi_{\Lambda_1}^{V \Lambda_0}(z) : 
 V(\Lambda_1) \longrightarrow V_{z}\otimes V(\Lambda_0),\\
&\Psi_{\Lambda_2}^{V \Lambda_1}(z) : 
 V(\Lambda_2) \longrightarrow V_{z}\otimes V(\Lambda_1),
\endaligned
\tag3.18 \\
&
\aligned
&\Phi_{\lambda_{\alpha}}^{\lambda_{\alpha+1} V^*}(z) :
 V(\lambda_{\alpha}) \longrightarrow V(\lambda_{\alpha+1})\otimes
 V_{z}^{*a} ,\\ 
&\Phi_{\Lambda_0}^{\Lambda_1 V^*}(z) :
 V(\Lambda_0) \longrightarrow V(\Lambda_1)\otimes
 V_{z}^{*a} ,\\ 
&\Phi_{\Lambda_1}^{\Lambda_2 V^*}(z) :
 V(\Lambda_1) \longrightarrow V(\Lambda_2)\otimes
 V_{z}^{*a} ,
\endaligned
\qquad
\aligned
&\Psi_{\lambda_{\alpha}}^{V^* \lambda_{\alpha+1}}(z) : 
 V(\lambda_{\alpha}) \longrightarrow V_{z}^{*a}\otimes 
 V(\lambda_{\alpha+1}), \\
&\Psi_{\Lambda_0}^{V^* \Lambda_1}(z) : 
 V(\Lambda_0) \longrightarrow V_{z}^{*a}\otimes 
 V(\Lambda_1), \\
&\Psi_{\Lambda_1}^{V^* \Lambda_2}(z) : 
 V(\Lambda_1) \longrightarrow V_{z}^{*a}\otimes 
 V(\Lambda_2), 
\endaligned \tag3.19
\endalign
$$
where $\lambda_{\alpha}=(1-\alpha)\Lambda_0+\alpha \Lambda_2$
for $\alpha \in {\Bbb R}$.
\endproclaim

\subhead Acknowledgments \endsubhead
The authors would like to thank H. Awata, M. Jimbo, A. Kuniba, 
T. Miwa, J. Suzuki, T. Takagi, A. Tsuchiya and Y. Yamada for 
stimulating discussions.
J. S. is very grateful to M. Wakimoto
for valuable discussions and kind hospitality 
while his visiting at Mie university in 1994.

\enddocument